\documentclass[12pt]{amsart}
\usepackage{amssymb,amsmath,amsthm}
\usepackage{epsfig}
\usepackage{graphicx}
\input xy
\xyoption{all}

\newtheorem{theorem}{Theorem}[section]
\newtheorem{lemma}[theorem]{Lemma}

\newtheorem{definition}[theorem]{Definition}

\newtheorem{corollary}[theorem]{Corollary}

\newtheorem{example}[theorem]{Example}

\long\def\symbolfootnote[#1]#2{\begingroup%
\def\thefootnote{\fnsymbol{footnote}}\footnote[#1]{#2}\endgroup}

\newcommand\Z{{\mathbb Z}}
\newcommand{\p}{\partial}
\newcommand{\sgn}{\text{sgn}}
\newcommand{\eop}{\hfill \ensuremath{\square} \\}
\begin{document}

\title{Entropic magmas, their homology, and related invariants of links and graphs}
\date{October 21, 2011}
\author{Maciej Niebrzydowski}
\address[Maciej Niebrzydowski]{University of Louisiana at Lafayette, Department of Mathematics,
	 217 Maxim D. Doucet Hall, 1403 Johnston Street,  Lafayette, LA 70504}
\email{mniebrz@gmail.com}
\author{J\'ozef H. Przytycki}
\address[J\'ozef H. Przytycki]{The George Washington University, Department of Mathematics,
	 Monroe Hall, Room 240, 2115 G Street NW,
	 Washington, D.C. 20052, and University of Gda\'nsk, Institute of Mathematics, Wita Stwosza 57, 80-952 Gda\'nsk, Poland} 
\email{przytyck@gwu.edu}

\thispagestyle{empty}

\begin{abstract} We define link and graph invariants from entropic magmas modeling them on the Kauffman bracket and Tutte polynomial. We define the homology of entropic magmas. We also consider groups that can be assigned to the families of compatible entropic magmas.
\end{abstract}

\maketitle

\section{Introduction}
A binary operation satisfying the entropic property $(a*b)*(c*d)=(a*c)*(b*d)$ was probably first considered by Sushkevich in 1937 \cite{Sus}.
Soon after Murdoch \cite{Mur} and Toyoda in a series of papers \cite{Toy-1,Toy-2,Toy-3,Toy-4} established the main properties of such magmas.
In particular, they proved the following result named after them; it is very important in our considerations.

\begin{theorem}\cite{Mur,Toy-2}
If $(X;*)$ is an entropic quasigroup, then $X$ has an abelian group structure such that $a*b=f(a)+g(b)+c$, where $f,g\colon X\to X$ are commuting group automorphisms.
\end{theorem}

The phrase {\it entropic property} that we use was coined in 1949 by Etherington \cite{Eth}. Other names for this property include: mediality,
bi-commutativity, alternation, bisymmetry, and abelianity. The word entropic refers to inner turning.

In the first part of the paper we show how to use entropic property for defining invariants of links and graphs that are based on the Kauffman bracket and Tutte polynomial so that the assumption of linearity is not present in their construction. The second part introduces homology of entropic magmas and its counterpart for families of entropic magmas connected by a condition $(a*_ib)*_j(c*_id)=(a*_jc)*_i(b*_jd)$.

\section{Tait graphs}\label{Tait}

We recall here (see e.g. \cite{Prz-37}) the construction of Tait which gives a bijection between signed plane 
graphs and link diagrams; bijection which also sends a Kauffman bracket version of the Tutte 
polynomial to the Kauffman bracket of a link diagram.

\begin{figure}
\begin{center}
\includegraphics[height=3cm]{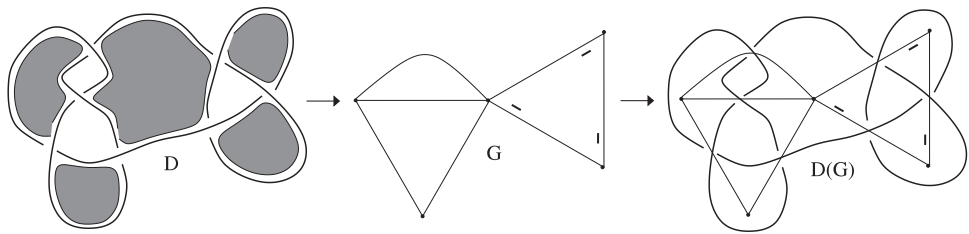}
\caption{Tait's construction of a graph from a link diagram (connected sum of the figure eight knot and 
the left handed trefoil knot with white infinite region), and back to the signed graph\label{kgk}}
\end{center}
\end{figure}

Tait was the first to notice the relation between knots and
planar graphs. He colored the regions of the knot diagram alternately white and black (motivated by Listing \cite{Lis}) and constructed the graph by placing a vertex inside each white region, and then connecting vertices by edges going through the crossing points of the diagram (we place however a vertex in black regions in our considerations, and study 
a ``black" graph instead of a dual ``white" graph of Tait; Figure \ref{kgk}).

\begin{figure}
\begin{center}
\includegraphics[height=2cm]{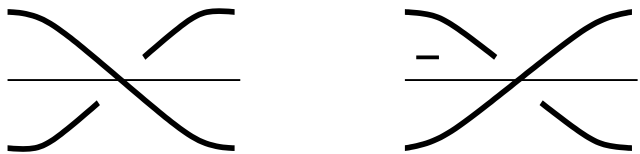}
\caption{convention for crossings assigned to signed edges (edges without markers are assumed to be positive)
\label{crossedges}}
\end{center}
\end{figure}

It is useful to mention the Tait construction going in the opposite direction, from
a signed planar graph, $G$, to a link diagram $D(G)$. We replace every edge of a graph by a crossing
according to the convention of Figure \ref{crossedges} and connect endpoints along edges as in
Figure \ref{kgk}.

\section{Entropic magma}\label{Entropic magma}
We define link diagram invariants from entropic magmas modeling them on the Kauffman bracket polynomial (but not assuming linearity). In 1984 the entropic magmas were used for the first time in Knot Theory; 
they were used to generalize 
Jones and Homflypt link invariants using Conway's skein triple \cite{P-T-1,P-T-2,Prz-1,Prz-100}.
If we generalize the Kauffman bracket relation, then in order to have link invariants (for framed links) 
from a magma $(A;*)$, we need an entropic condition,
and an additional condition for the second Reidemeister move. We analyze it below and if initial conditions
are given by $(a_1,\ldots,a_n,\ldots)$, we call it a bracket (or Kauffman bracket) magma. This we discuss below 
(later we consider signed graphs and play analogous game with invariants coming from Tutte magma).

\subsection{Kauffman bracket magma}\label{kauffbrac}
\begin{figure}
\begin{center}
\includegraphics[height=3.5cm]{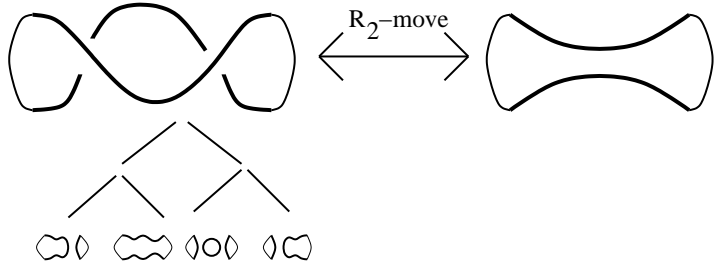}
\caption{$R_2$-move; denominator closure
\label{R2-denom}}
\end{center}
\end{figure}

\begin{figure}
\begin{center}
\includegraphics[height=5cm]{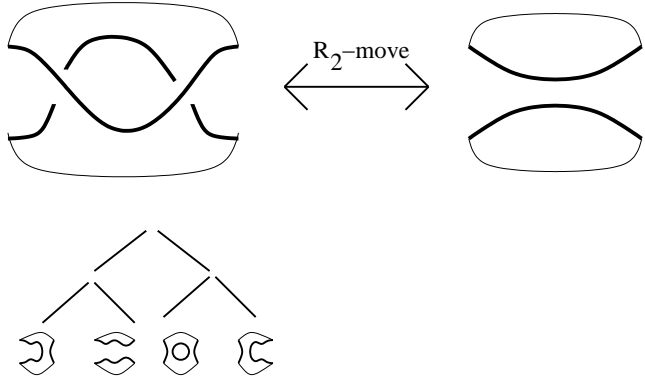}
\caption{$R_2$-move; numerator closure
\label{R2-numer}}
\end{center}
\end{figure}

We consider link invariants taking values in an entropic magma $(A;*)$ with a chosen sequence of not necessarily different elements $a_1,a_2,\ldots,a_n,\ldots$, which satisfy the following conditions for any $n\geq 1$. The relations are derived from the second Reidemeister 
moves (assuming that the crossings are ordered and the two new crossings are the last). We get (i) or (ii) as below:\\
$\mbox{ (i) } (a_{n+1}*a_n)*(a_{n+2}*a_{n+1}) = (a_{n+1}*a_{n+2})*(a_n*a_{n+1}) = a_n;$\\ 
this case, for $n=1$, is illustrated in Figure \ref{R2-denom},\\ 
$\mbox{ (ii) } (a_{n}*a_{n+1})*(a_{n+1}*a_{n}) = a_{n+1};$\\
this case, for $n=1$, is illustrated in Figure \ref{R2-numer}.

First, we construct a regular isotopy\footnote{In regular isotopy we consider unoriented 
link diagrams modulo the second and the third Reidemeister moves. In an equivalent approach, we 
can consider unoriented framed links and use the fact that two diagrams on $S^2$ representing 
framed links with blackboard framing are ambient isotopic iff they are related by 
$R_2$ and $R_3$ moves. Notice that $R_1$ is changing framing by $\pm 1$.} link invariant from a Kauffman 
bracket magma:

\begin{theorem}
Let $(A;*,a_1,a_2,\ldots)$ be a Kauffman bracket magma.
Then there is a unique unoriented framed link invariant, say $P: \{Links^{fr}\} \to A$, such that 
$P(T_n)=a_n$ and if $L=$ \parbox{0.9cm}{\psfig{figure=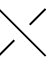,height=0.7cm}}, 
$L_0=$ \parbox{0.9cm}{\psfig{figure=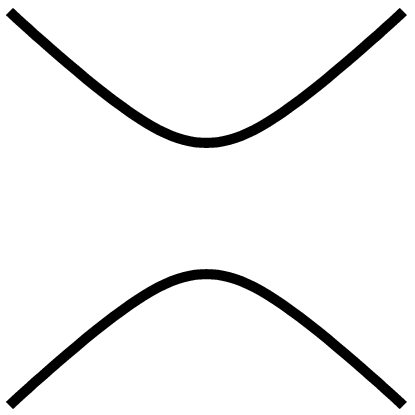,height=0.7cm}}, and 
$L_{\infty}=$ \parbox{0.9cm}{\psfig{figure=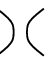,height=0.7cm}} is a Kauffman bracket skein triple, then 
$$ P(L)=P(L_0)*P(L_{\infty}).$$ 
\end{theorem}
\begin{proof} Our proof follows the Kauffman's proof of existence of the Kauffman bracket polynomial but 
we give it here for completeness. The proof is organized as follows:\\
(1) We consider any link diagram $D$ representing $L$ and we order its crossings. 
Then, we compute (in a unique way) the value (in $A$) of the diagram, resolving crossings according 
to their ordering and using initial values at the end (at leaves of the computational tree).
We denote this invariant of diagrams with ordered crossings by $P(D)$.\\
(2) We show that obtained value does not depend on the ordering of crossings. 
Here, it suffices to check that we can exchange the order of two crossings which were consecutive 
in the initial ordering, say $v_1<v_2$. 
For a diagram $D$ with crossings $v_1,v_2$, we write $D^{v_1,v_2}_{\varepsilon_1,\varepsilon_2}$
for a diagram obtained from $D$ by smoothing $v_1$ and $v_2$ according to $\varepsilon_1,\varepsilon_2$,
where $\varepsilon_i=0 \mbox{ or } \infty$ decides which Kauffman marker we use for smoothing; 
of course $D^{v_1,v_2}_{\varepsilon_1,\varepsilon_2}= D^{v_2,v_1}_{\varepsilon_2,\varepsilon_1}$.
We have (in an initial order):
$$P(D)= P(D^{v_1}_{0})*P(D^{v_1}_{\infty})=$$
$$(P(D^{v_1,v_2}_{0,0})*P(D^{v_1,v_2}_{0,\infty}))*(P(D^{v_1,v_2}_{\infty,0})*P(D^{v_1,v_2}_{\infty,\infty}))$$
After switching the order of $v_1$ and $v_2$ we smooth $v_2$ first and obtain analogously:
$$P(D)= P(D^{v_2}_{0})*P(D^{v_2}_{\infty})=$$
$$(P(D^{v_2,v_1}_{0,0})*P(D^{v_2,v_1}_{0,\infty}))*(P(D^{v_2,v_1}_{\infty,0})*P(D^{v_2,v_1}_{\infty,\infty}))=$$
$$(P(D^{v_1,v_2}_{0,0})*P(D^{v_1,v_2}_{\infty, 0}))*(P(D^{v_1,v_2}_{0,\infty})*P(D^{v_1,v_2}_{\infty,\infty}))$$
Our invariants coincide on four diagrams $D^{v_1,v_2}_{0,0},D^{v_1,v_2}_{0,\infty}, D^{v_1,v_2}_{\infty,0}$ and 
$D^{v_1,v_2}_{\infty,\infty}$. Thus, in order for them to coincide on $D$ it suffices to have 
entropic condition, in our case:
$$(P(D^{v_1,v_2}_{0,0})*P(D^{v_1,v_2}_{0,\infty}))*(P(D^{v_1,v_2}_{\infty,0})*P(D^{v_1,v_2}_{\infty,\infty}))=$$ 
$$(P(D^{v_1,v_2}_{0,0})*P(D^{v_1,v_2}_{\infty, 0}))*(P(D^{v_1,v_2}_{0,\infty})*P(D^{v_1,v_2}_{\infty,\infty}))$$
is a special case of the entropic condition.\\
(3) We analyze the behavior of our diagram invariant under the second Reidemeister move $R_2$. Because (by 2) we can 
put crossings involved in $R_2$ at the very end of the calculations, we conclude that the relations (i) and (ii) 
of Figure \ref{R2-denom} and \ref{R2-numer} are necessary and sufficient.\\
(4) The invariance under the third Reidemeister move follows from (3), in a similar manner like 
in the case of the Kauffman bracket -- we start the resolution from the top crossing; we illustrate it in 
Figure \ref{R3Kauff}.

\begin{figure}
\begin{center}
\includegraphics[height=4.7cm]{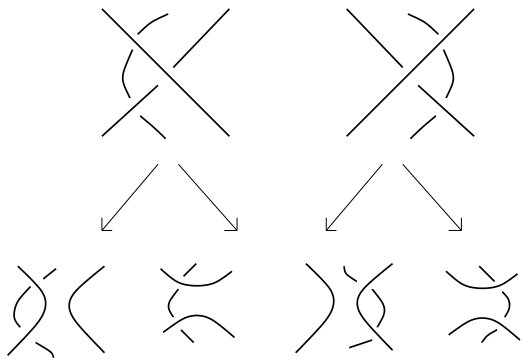}
\caption{After smoothings of upper crossings, diagrams differ by $R_2$ moves only.
\label{R3Kauff}}
\end{center}
\end{figure}

More precisely, by (2), we can start our calculation from the upper crossing involved in the third Reidemeister move.
After the first resolution, the left diagrams differ by two $R_2$ moves (so we can use (3)), and the 
right diagrams are isotopic (Figure \ref{R3Kauff}).
\end{proof}

\subsection{Kauffman bracket magma invariants and the first Reidemeister moves}

We discuss in this section the change of $P(D)$ under first Reidemeister moves.
We have to consider two types of the first move, a positive $R_{+1}$ in which we add a 
positive kink, and a negative $R_{-1}$ in which a negative kink is added to a diagram $D$.
Equivalently, we measure the effect on $P(L)$ of twisting a framing in a positive 
direction ($L \to L^{(1)}$) or the negative one ($L \to L^{(-1)}$). Denote by $\tau(D)$ 
a link diagram $D \sqcup \bigcirc$, that is, we add a trivial component to the diagram $D$.
Then we have $P(D^{(1)})=P(\includegraphics[trim=0mm 0mm 0mm 0mm, width=.02\linewidth]
{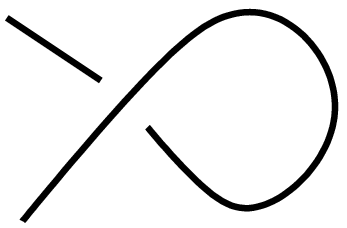}) = P(\tau(D))*P(D)$, and analogously:
$P(D^{(-1)})=P(\includegraphics[trim=0mm 0mm 0mm 0mm, width=.02\linewidth]
{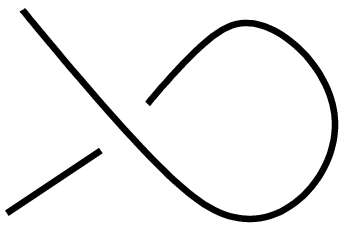}) = P(D)*P(\tau(D))$.

For example, a positive loop (i.e. the trivial knot with framing $1$, $T_1^{(1)}$= 
\parbox{1.4cm}{\psfig{figure=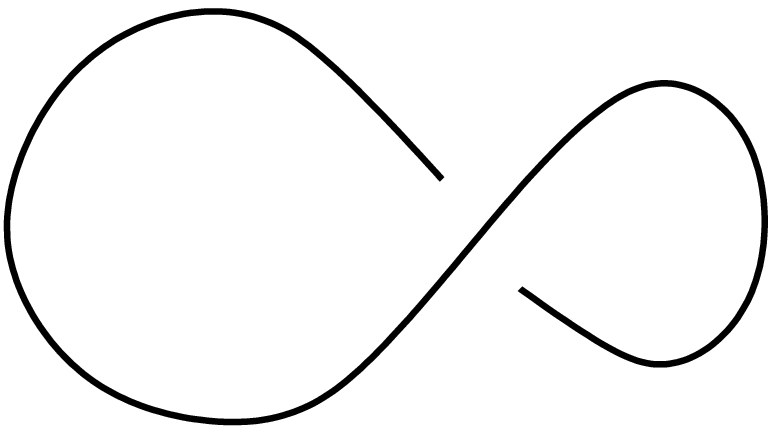,height=0.7cm}}) 
has value $a_2*a_1$, while $T_1^{(-1)}=$ \parbox{1.4cm}{\psfig{figure=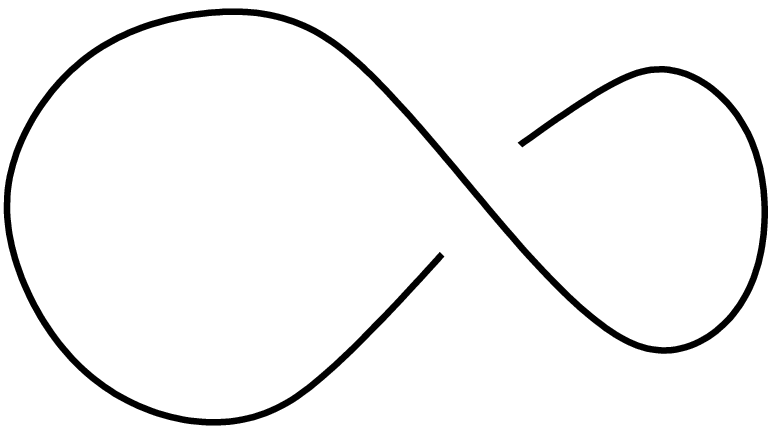,height=0.7cm}}
 has the value $a_1*a_2$.

Notice that the equality $P((L^{(1)})^{(-1)})=P(L)$ follows from the relation (i) of Section \ref{kauffbrac} even if twisting is performed on different component of $L$. Namely,
we can assume that kink crossings are the last crossings of the computational tree, thus we can 
assume that $L$ is represented by a trivial diagram of $n$ components.
Then we have $P(L^{(1)})= a_{n+1}*a_n$,  $P(L^{(-1)})= a_{n}*a_{n+1}$, and 
$P((L^{(1)})^{(-1)})= (a_{n+1}*a_{n+2})*(a_n*a_{n+1}) \stackrel{(i)}{=}  a_n=P(L)$.

Let us extend the notation $\tau(L)=L\sqcup \bigcirc$ to a subset of $A$ realized by link invariants.
More precisely, let $A^{\mathcal L}= \{a\in A \ | \ a=P(D)\ \mbox{for some diagram D }\}$.
Then on $A^{\mathcal L}$ we define $\tau(a) = P(\tau(D))$ where $P(D)=a$.
It is not always the case that $\tau$ can be extended from $A^{\mathcal L}$ to $A$, but if 
it does extend, we still denote an extension by $\tau$.

To produce invariant of unoriented links (diagrams under all Reidemeister moves), we can assume that we 
consider only zero framings on links (that is, the framing given by Seifert surfaces for each component),
or that we orient $L$ and consider framing given by a Seifert surface of oriented $L$.
This can be rephrased in the language using Tait or writhe number of an oriented diagram 
(the sum of signs of crossings), as was done in the case of Kauffman bracket (modified to 
Jones polynomial). 

\begin{figure}
\begin{center}
\includegraphics[height=4.3cm]{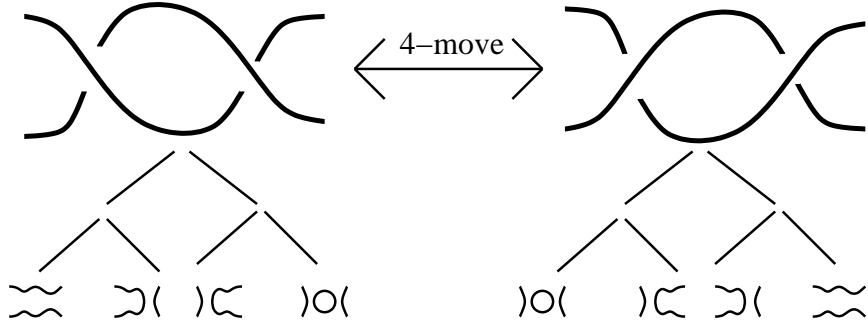}
\caption{4-move and resolving computational trees
\label{4-moveAlg}}
\end{center}
\end{figure}

\begin{lemma}
Consider a Kauffman bracket magma $(A;*,a_1,a_2,\ldots)$ which also satisfies, for every $n\geq 1$:
$$ (a_{n}*a_{n+1})*(a_{n+1}*a_{n+2})= (a_{n+2}*a_{n+1})*(a_{n+1}*a_{n})$$
Then the invariant of links $P(L)$ is preserved by 4-moves.
\end{lemma}
\begin{proof} We illustrate in Figure \ref{4-moveAlg} how the 4-move is changing the Kauffman bracket magma invariant.
This involves two variables and can be read as
$(a*b)*(b*\tau(b))= (\tau(b)*b)*(b*a)$. But in fact, as we work with diagrams on $R^2$ we can smooth all 
crossings not involved in the move first and then deal exclusively with a clasp, again having 
two cases, numerator and denominator, as illustrated in Figure \ref{4-moveAlg}.
\end{proof}

We should also observe that if we work modulo relations giving invariants of 4-moves, then 
the invariant is preserved by framing change of multiplicity of four.

\section{Tutte magma for signed graphs}
We are motivated by a bijection between signed plane graphs and link diagrams (on $R^2$ or on $S^2$, whatever is 
more convenient), but we define magma invariants for all signed graphs, for any magma $(A;*)$ with a chosen sequence of elements
$a_1,a_2,\ldots,a_n,\ldots$ (not necessarily all different).

If edges of $G$ are ordered, we define an invariant of a signed graph in Tutte fashion, 
by giving to a graph with no edges
and $n$ vertices (say $T_n$) the value $a_n$ (that is $P(T_n)=a_n$), and then:\\
(i) For a positive  edge $e_+$ which is not a loop:
$$ P(G)=P(G/e_+)*P(G-e_+),$$
here, $G-e_+$ and $G/e_+$ denote deleting and contracting the edge $e_+$, respectively.\\
(ii) For a negative edge $e_-$ which is not a loop:
$$ P(G)= P(G-e_-)* P(G/e_-)$$
(iii) If $e_+$ is a positive loop then:
$$ P(G)=P(G//e_+)*P(G-e_+),$$$$\mbox{ where $G//e_+$ is obtained from $G/e_+$ by adding an isolated vertex};$$
(iv) If $e_-$ is a negative loop then:
$$ P(G)= P(G-e_-)* P(G//e_-),$$$$\mbox{ where $G//e_-$ is obtained from $G/e_-$ by adding an isolated vertex}.$$

When computing the value of $P(G)\in A$, we use the edges of $G$ one by one according to their ordering.

\begin{figure}
\begin{center}
\includegraphics[height=3.7cm]{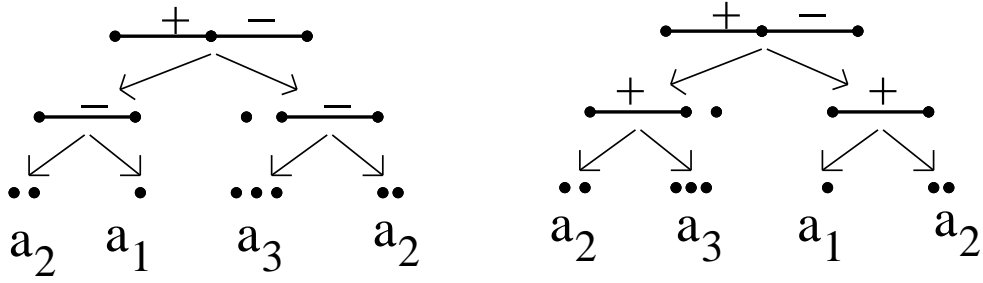}
\caption{changing the ordering of edges in a computational tree
\label{Gpl-min-tree}}
\end{center}
\end{figure}

The invariant is well defined but depends on the ordering of edges. 
If we switch the order of two edges, we may have different result: one of the form $(a*b)*(c*d)$, 
and the second of the form $(a*c)*(b*d)$. Thus, entropic magma may be used to get an invariant of 
a signed graph. We can try to do better and consider only some entropic relations.
For example, if $G$ has 2 edges, the only type of relation which $(A;*)$ has to satisfy is an entropic relation:
$$(a_{n+1}*a_n)*(a_{n+2}*a_{n+1}) = (a_{n+1}*a_{n+2})*(a_n*a_{n+1}).$$
We meet this condition when considering the graph \parbox{2.9cm}{\psfig{figure=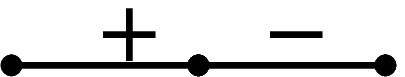,height=0.4cm}} 
with additional $n-1$ isolated points. Its two computational trees are shown in Figure \ref{Gpl-min-tree}.

When analyzing graphs with 3 edges, we find 4 additional entropic conditions:\\
By considering the graph \ \parbox{3.9cm}{\psfig{figure=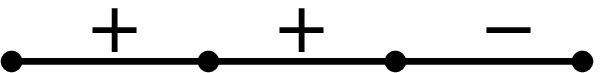,height=0.4cm}} we get
$$((a_{n+2}*a_{n+1})*(a_{n+3}*a_{n+2}))*((a_{n+1}*a_{n})*(a_{n+2}*a_{n+1}))=$$
$$((a_{n+2}*a_{n+1})*(a_{n+1}*a_{n}))*((a_{n+3}*a_{n+2})*(a_{n+2}*a_{n+1}))$$ 
and its reverse (by considering \ \parbox{3.9cm}{\psfig{figure=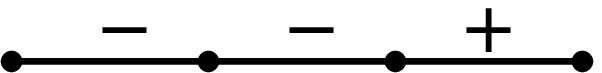,height=0.4cm}}):
$$ ((a_{n+1}*a_{n+2})*(a_{n}*a_{n+1}))*((a_{n+2}*a_{n+3})*(a_{n+1}*a_{n+2}))=$$
$$ ((a_{n+1}*a_{n+2})*(a_{n+2}*a_{n+3}))*((a_{n}*a_{n+1})*(a_{n+1}*a_{n+2})).$$
From the positive graph \parbox{2.0cm}{\psfig{figure=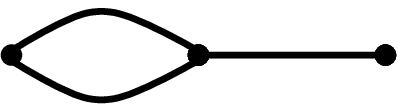,height=0.5cm}} 
we get an entropic identity:
$$ ((a_{n+1}*a_{n})*(a_{n}*a_{n+1}))*((a_{n+2}*a_{n+1})*(a_{n+1}*a_{n+2}))=$$
$$ ((a_{n+1}*a_{n})*(a_{n+2}*a_{n+1}))*((a_{n}*a_{n+1})*(a_{n+1}*a_{n+2}))$$
and its reverse (by considering the negative graph \parbox{2.0cm}{\psfig{figure=Gplpl-pl.eps,height=0.5cm}}):
$$ ((a_{n+2}*a_{n+1})*(a_{n+1}*a_{n+2}))*((a_{n+1}*a_{n})*(a_{n}*a_{n+1}))=$$
$$ ((a_{n+2}*a_{n+1})*(a_{n+1}*a_{n}))*((a_{n+1}*a_{n+2})*(a_{n}*a_{n+1})).$$

In general, we have to consider all graphs, so there are no easy criteria for setting the conditions. Thus, we define Tutte magma abstractly 
as a magma with distinguished elements $a_1,a_2,\ldots,a_n$, such that if we compute invariants of any fixed signed graph 
$G$ using different orders of the edges, we always get the same result (we know that we need some, but not all entropic conditions).  
Formally:
\begin{definition} Let ${\mathcal A}=(A,*,a_1,a_2,\ldots)$ be a magma with a chosen sequence of elements,
and $P(G^{ord})$ an associated invariant of graphs with ordered edges. 
We say that ${\mathcal A}$ is a Tutte magma if $P(G^{ord})$ does not depend on the ordering of edges of $G$, 
for any signed graph $G$.
\end{definition}
Clearly, if ${\mathcal A}$ is an entropic magma, it is also a Tutte magma, but as we demonstrated earlier, 
not all entropic relations are needed. We can rephrase the definition using partial computational trees
for two edges with different order, but we leave it to the reader.
Similarly, we can define partial Tutte magma (by analogy to partial Conway algebra in \cite{P-T-2,Prz-100}).

Now we describe another invariant of signed graphs, using the function $P$ and congruences on algebras. First, we recall the notion of a congruence on a magma.
\begin{definition} Let ${\mathcal A}=(A,*,a_1,a_2,\ldots)$ be a magma, and let $\theta$ be an equivalence relation on $A$. $\theta$ is a congruence on ${\mathcal A}$ if $b_1\theta c_1$ and $b_2\theta c_2$ implies $(b_1*b_2)\theta (c_1*c_2)$, for $b_1$, $b_2$, $c_1$, $c_2\in A$.
\end{definition}
\begin{definition} Given a magma ${\mathcal A}=(A,*,a_1,a_2,\ldots)$ and a congruence $\theta$ on ${\mathcal A}$, one can form a quotient magma
${\mathcal A}/\theta$, whose underlying set is the set of equivalence classes of $\theta$, with the binary operation
$[a]*[b]=[a*b]$, where $[x]$ denotes the equivalence class of $x\in A$.
\end{definition}

Now we are ready to define the invariant.

\begin{definition}
We begin by fixing a magma ${\mathcal A}=(A,*)$ with a sequence of not necessarily different elements $a_1,a_2,\ldots$. For any signed graph $G$, let $X_G\subseteq A$ denote the set of values of the function $P$ obtained using all possible orderings of the edges of $G$. Let $\theta$ be the smallest congruence on ${\mathcal A}$ generated by the set $X_G$ (that is, the smallest congruence containing all the pairs $(a,b)$ for $a$, $b\in X_G$). Form the quotient algebra 
${\mathcal A}_G={\mathcal A}/\theta$. Then the function $T$ assigning a quotient algebra ${\mathcal A}_G$ to a signed graph $G$ is an invariant of signed graphs.
\end{definition}

\subsection{Examples of calculations for line graphs and polygons}

This section contains calculations of $P$ for certain signed plane graphs. We keep in mind the correspondence between such graphs and link diagrams described in Section \ref{Tait} (especially when using $\tau$ below).

\begin{example} Let $L_n$ be a line graph (that is, a tree with $n+1$ vertices, all except two of degree two), with positive edges.
Then we have $$P(L_n)= P(L_{n-1})* \tau(P(L_{n-1}))$$ with an initial condition: $P(L_0)=a_1$. 
We have for example:
$$P(L_1)= P(L_0)*\tau(P(L_0))= a_1*a_2$$
$$P(L_2)= P(L_1)*\tau(P(L_1))= (a_1*a_2)*(a_2*a_3)$$
$$P(L_3)= P(L_2)*\tau(P(L_2))= ((a_1*a_2)*(a_2*a_3))*((a_2*a_3)*(a_3*a_4))$$
$$P(L_4)= P(L_3)*\tau(P(L_3))= $$ 
$$(((a_1*a_2)*(a_2*a_3))*((a_2*a_3)*(a_3*a_4)))*(((a_2*a_3)*(a_3*a_4))*((a_3*a_4)*(a_4*a_5)))$$
\end{example}

\begin{example}
Let $C_n$ be an $n$-gon (a cycle graph with $n$ edges and $n$ vertices). Then we have 
$$P(C_n)= P(C_{n-1})* P(L_{n-1})$$ with an initial condition: $P(C_1)=a_2*a_1$.
We have for example:
$$P(C_2)=P(C_1)*P(L_1)= (a_2*a_1)*(a_1*a_2)$$
$$P(C_3)=P(C_2)*P(L_2)= ((a_2*a_1)*(a_1*a_2))*((a_1*a_2)*(a_2*a_3))$$
$C_3$ corresponds to a right handed trefoil knot ( $\bar 3_1$ =  \parbox{0.9cm}{\psfig{figure=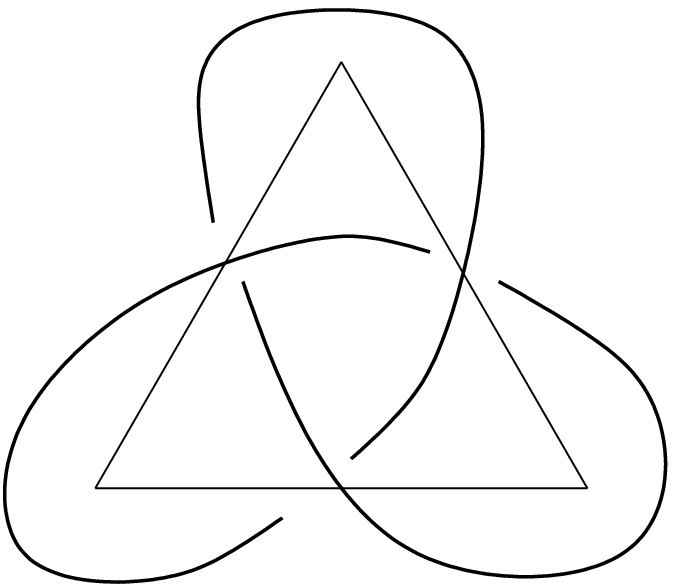,height=0.7cm}}).\\
The notation $3_1$ is used in the Rolfsen's table of knots for the left handed trefoil knot 
\parbox{0.9cm}{\psfig{figure=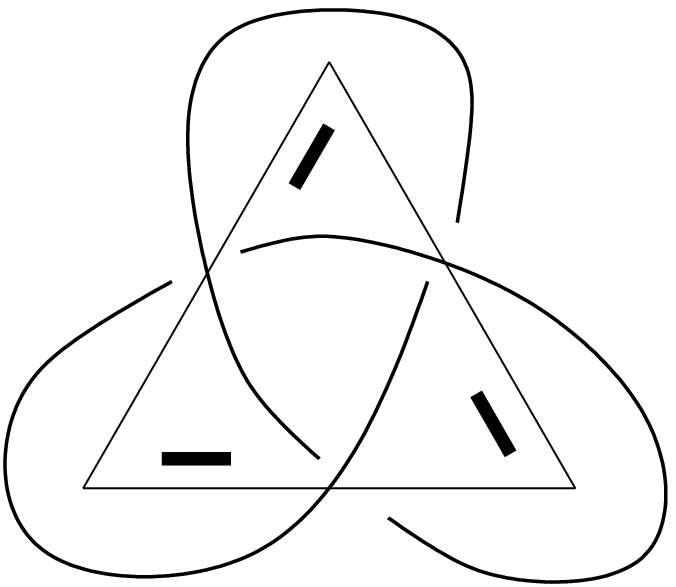,height=0.7cm}}.
\end{example}
\begin{example}
Let $C'_n$ be an $n$-gon with one edge doubled, then we have:
$$P(C'_n)= (\tau(P(C_{n-1}))*P(C_{n-1}))*P(C_n)$$
For example $C'_2$ is a theta curve and describes the left handed trefoil knot; we get:
$$P(3_1)=P(C'_2)= (\tau(P(C_{1}))*P(C_{1}))*P(C_2)= $$$$((a_3*a_2)*(a_2*a_1))*((a_2*a_1)*(a_1*a_2))$$
The graph $C'_3$ corresponds to the figure eight knot \parbox{0.9cm}{\psfig{figure=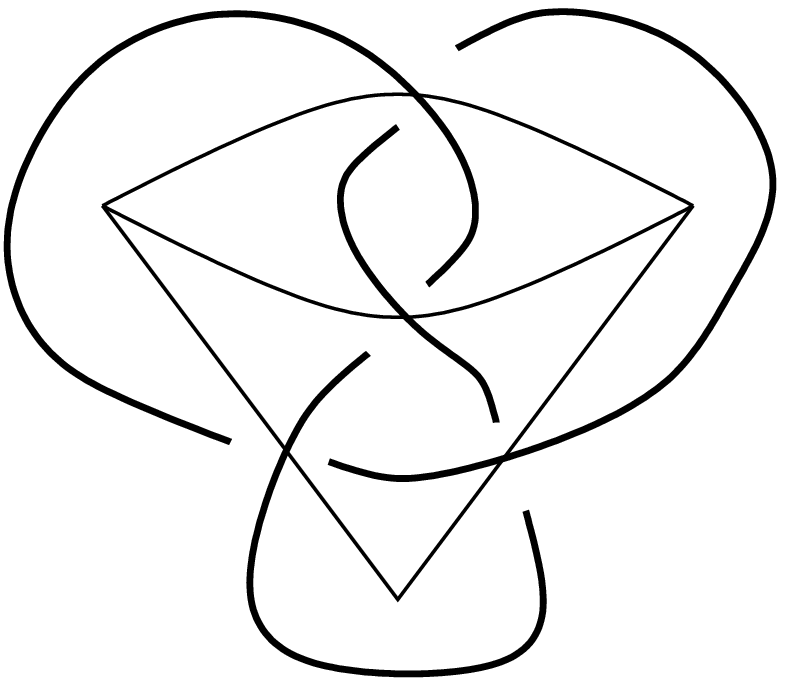,height=0.7cm}},
and we have:
$$P(4_1)= P(C'_3)= (\tau(P(C_{2}))*P(C_{2}))*P(C_3)= $$
$$(((a_3*a_2)*(a_2*a_3))*((a_2*a_1)*(a_1*a_2)))*(((a_2*a_1)*(a_1*a_2))*((a_1*a_2)*(a_2*a_3))).$$
\end{example}

\section{Homology of entropic magmas}
Let $(X;*)$ be a magma, $A$ any set and $\pi\colon A\times X\to X$ the projection to
the second coordinate. Any magma structure on $A\times X$ for which $\pi$ is an epimorphism
can be given by a system of functions $\phi_{a_1,a_2}(x_1,x_2)\colon X\times X\to A$ by:
$$(a_1,x_1)*(a_2,x_2)= (\phi_{a_1,a_2}(x_1,x_2), x_1*x_2).$$
Functions $\phi_{a_1,a_2}(x_1,x_2)$ are uniquely defined by the multiplication on $A\times X $.
Thus, binary operations on $A\times X$ agreeing with $\pi$ are in bijection with choices of functions
$\phi_{a_1,a_2}$.
If we require some special structure on  $(X;*)$ (e.g. associativity or entropic condition), we obtain
some property of $\phi_{a_1,a_2}(x_1,x_2)$ which we call a dynamical cocycle for the structure.
\begin{enumerate}
\item[(1)] Let $(X;*)$ be a semigroup; in order that an action on $A\times X$ be associative, we
need:
$$((a_1,x_1)*(a_2,x_2))*(a_3,x_3)=(\phi_{a_1,a_2}(x_1,x_2), x_1*x_2)*(a_3,x_3)=$$
$$ (\phi_{\phi_{a_1,a_2}(x_1,x_2),a_3}(x_1*x_2,x_3), (x_1*x_2)*x_3)$$
to be equal to
$$(a_1,x_1)*((a_2,x_2)*(a_3,x_3))= (a_1,x_1)*(\phi_{a_2,a_3}(x_2,x_3), x_2*x_3)= $$
$$ (\phi_{a_1,\phi_{a_2,a_3}(x_2,x_3)}(x_1,x_2*x_3),x_1*(x_2*x_3)).$$
Thus, the dynamical cocycle condition in the associative case has a form:
$$\phi_{\phi_{a_1,a_2}(x_1,x_2),a_3}(x_1*x_2,x_3)= \phi_{a_1,\phi_{a_2,a_3}(x_2,x_3)}(x_1,x_2*x_3). $$
\item[(2)] We now assume that $(X;*)$ is an entropic magma, that is
$(a*b)*(c*d)= (a*c)*(b*d)$ for any $a,b,c,d \in X$. We look for the condition on dynamical cocycle
so that $A\times X$ is entropic. We need
$$((a_1,x_1)*(a_2,x_2))*((a_3,x_3)*(a_4,x_4))=$$
$$  (\phi_{a_1,a_2}(x_1,x_2), x_1*x_2)* (\phi_{a_3,a_4}(x_3,x_4), x_3*x_4)=$$
$$ (\phi_{\phi_{a_1,a_2}(x_1,x_2),\phi_{a_3,a_4}(x_3,x_4)}( x_1*x_2,x_3*x_4), (x_1*x_2)*(x_3*x_4))$$
to be equal to
$$((a_1,x_1)*(a_3,x_3))*((a_2,x_2)*(a_4,x_4))=$$
$$  (\phi_{a_1,a_3}(x_1,x_3), x_1*x_3)* (\phi_{a_2,a_4}(x_2,x_4), x_2*x_4)=$$
$$ (\phi_{\phi_{a_1,a_3}(x_1,x_3),\phi_{a_2,a_4}(x_2,x_4)}( x_1*x_3,x_2*x_4), (x_1*x_3)*(x_2*x_4)).$$
Thus, the dynamical cocycle condition in entropic case has the form:
$$\phi_{\phi_{a_1,a_2}(x_1,x_2),\phi_{a_3,a_4}(x_3,x_4)}( x_1*x_2,x_3*x_4)=
  \phi_{\phi_{a_1,a_3}(x_1,x_3),\phi_{a_2,a_4}(x_2,x_4)}( x_1*x_3,x_2*x_4).$$
\end{enumerate}
Extensions of modules, groups and Lie algebras are described in the
classical book by Cartan and Eilenberg \cite{C-Eil}.

Let $(X;*)$ be an entropic magma, $A$ an abelian group with a given pair of commuting
homomorphisms $t,s\colon A\to A$ and a constant $a_0\in A$; we consider $(A;*)$ as an entropic magma
with an affine action $a*b= ta +sb + a_0$.
Then we define a binary operation on $A\times X$ by
$(a_1,x_1)*(a_2,x_2) = (a_1*a_2 + f(x_1,x_2),x_1*x_2)$.
In order for $A\times X$ to be an entropic magma,
$\phi_{a_1,a_2}(x_1,x_2) = a_1*a_2 + f(x_1,x_2)$ should be an entropic dynamical cocycle.
This leads to entropic cocycle condition:
$$tf(x_1,x_2) - tf(x_1,x_3) +sf(x_3,x_4) - sf(x_2,x_4) +$$$$
f(x_1*x_2,x_3*x_4) - f(x_1*x_3,x_2*x_4)=0.$$
The above formula served as a hint in defining entropic homology.
In particular, for a unital ring $R$, $\partial\colon RX^4\to RX^2$ may be given by:
$$\partial(x_1,x_2,x_3,x_4)= 
t(x_1,x_2) - t(x_1,x_3) + s(x_3,x_4) - s(x_2,x_4) +$$$$ (x_1*x_2,x_3*x_4) - (x_1*x_3,x_2*x_4)$$
and $\partial\colon RX^2 \to RX$  may be given by:
$\partial(x_1,x_2)= tx_1 - x_1*x_2 + sx_2$. 

The last map (on the level of cohomology) is derived as follows (we follow the classical 
case of a group extension by an abelian group; we adjust it to an extension of an entropic magma $(X;*)$  
by an entropic, affine magma $(A;*)$ ):\\
We ask when extensions $A\times X \to X$ given by various $f\colon X\times X \to A$ with
$(a_1,x_1)*(a_2,x_2) = (a_1*a_2 + f(x_1,x_2),x_1*x_2)$ are equivalent in the ``fiber preserving" sense.
That is, we would like to know whether two extensions $A\times X \to X$ with the same affine action on $A$ ($a*b= ta +sb + a_0$) 
are related by a magma homomorphism $F\colon A\times X \to A\times X$, where $F$ is constant on $X$ factor.
We can express $F$ as $F(a,x)= (c(x)+a,x)$, with $c\colon X \to A$. Homomorphism means, of course, that:
$$F((a_1,x_1)*_1(a_2,x_2))= F(a_1,x_1) *_2  F(a_2,x_2).$$ Here $*_1$ corresponds to $f_1$ and $*_2$ to $f_2$,
so we write concretely:
$$ (a_1*a_2+ f_1(x_1,x_2)+ c(x_1*x_2),x_1*x_2)=  (a_1+ c(x_1),x_1)*_2(a_2+ c(x_2),x_2)$$
and further:
$$(a_1*a_2+ f_1(x_1,x_2)+ c(x_1*x_2),x_1*x_2)=$$$$ ((a_1+ c(x_1))*(a_2+ c(x_2)) + f_2(x_1,x_2), x_1*x_2)=$$
$$(t(a_1+ c(x_1))+ s(a_2+ c(x_2)) + a_0 + f_2(x_1,x_2), x_1*x_2)= $$$$(a_1*a_2+ tc(x_1)+sc(x_2) + f_2(x_1,x_2), x_1*x_2).$$
Equivalently:
$$(f_1-f_2)(x_1,x_2)= tc(x_1)+sc(x_2) - c(x_1*x_2).$$ This suggests 
$(\partial c) (x_1,x_2)= tc(x_1)+sc(x_2) - c(x_1*x_2)$ and $(f_1-f_2)= \partial c$.
This further suggests $\partial\colon RX^2\to RX$ given by:
$$\partial(x_1,x_2)= tx_1 - x_1*x_2 + sx_2,$$ as we wrote before.
We should mention here that our map $\partial\colon Hom(RX^1,A) \to Hom(RX^2,A)$ given by
$$(\partial c) (x_1,x_2)= tc(x_1)+sc(x_2) - c(x_1*x_2)$$ 
was not using any specific properties of the magma $(X;*)$, thus for $s=1-t$ it gives the
coboundary operator for the twisted rack cohomology (see \cite{C-E-S}):
$$(\partial c) (x_1,x_2)= tc(x_1)+(1-t)c(x_2) - c(x_1*x_2)$$
and in the group case (but with the trivial group action giving a central extension), we put $t=s=1$ to get:
$$(\partial c) (x_1,x_2)= c(x_2) - c(x_1*x_2) +  c(x_1).$$

Now we will give a general definition of entropic homology which extends the above for s=t=0.
For a given entropic magma $(X;*)$, and $n\geq 0$, let $C_n(X)$ be the free abelian 
group generated by $2^n$-tuples $(x_1,x_2,\ldots,x_{2^n})$ of elements of $X$; 
in other words, $C_n(X) = {\Z}X^{2^n} = ({\Z}X)^{\otimes 2^n}$. We set $C_n(X)=0$ for $n<0$.
We are going to define several boundary homomorphisms going from $C_n(X)$ to $C_{n-1}(X)$.

\begin{figure}
\begin{center}
\includegraphics[height=4.7cm]{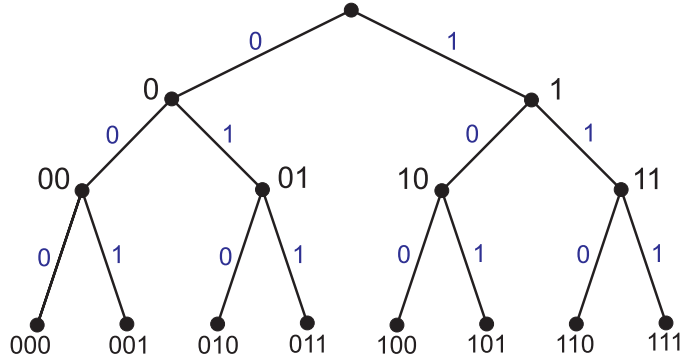}
\caption{Addresses of nodes in a binary tree\label{tree8}}
\end{center}
\end{figure}

Let $T_n$ be the $2^n$-tuple of addresses of nodes on the $n$-th level of the infinite binary tree. For example,
$T_3=\{000,001,010,011,100,101,110,111\}$ (see Figure \ref{tree8}). We view $T_n$ as a $2^n$-element set with elements numbered 1 through $2^n$. For brevity, we write $$w[i]_1=k$$ when the $i$-th place on the list $w$ has an address containing $k$ ones. Thus, for example, $T_3[4]_1=2$.
Let $S^k_{2^n}$ denote the set of permutations of the symmetric group $S_{2^n}$ that permute only addresses containing $k$ 1s. For example,
$S^1_{2^3}=\{(),(3,5),(2,3),(2,3,5),(2,5,3),(2,5)\}$, where $()$ is the identity permutation. For a $2^n$-tuple $w$ of elements of $X$,
let $\sigma w$ be the result of applying the permutation $\sigma$ to $w$; we extend this linearly to the whole $C_n$. We define a homomorphism $\delta^k_n\colon \Z X^{2^n}\to \Z X^{2^{n}}$ by $$\delta^k_n(w)=\sum_{\sigma \in S^k_{2^n}} \sgn(\sigma) \sigma w,$$
where $\sgn(\sigma)$ denotes the sign of the permutation $\sigma$.
Also, for a list $\nu$ of integers $(a_1,\ldots,a_{n-1})$,
we define $\delta^\nu_n\colon \Z X^{2^n}\to \Z X^{2^{n}}$ by $$\delta^\nu_n(w)=\sum_{i\in\{1,\ldots,n-1\}}a_i\delta^i_n(w).$$ For example, $\delta^{(2,1)}_3=3()+2(2,3,5)+2(2,5,3)-2(3,5)-2(2,3)-2(2,5)+(4,6,7)+(4,7,6)-(6,7)-(4,6)-(4,7)$.\\
We will also need a homomorphism $\mu_n\colon \Z X^{2^n}\to \Z X^{2^{n-1}}$ defined on $2^n$-tuples by:
$$\mu_n(x_1,x_2,x_3,x_4,\ldots,x_{2^n-1},x_{2^n})=(x_1*x_2,x_3*x_4,\ldots,x_{2^n-1}*x_{2^n}).$$
We will show that the maps $\p^k_n=\mu_n\delta^k_n$, with $k\in\{1,\ldots,n-1\}$, play the role of basic differentials.

\begin{theorem}\label{mainhom}
For $k\in\{1,\ldots,n-1\}$, $l\in\{1,\ldots,n-2\}$ and $n>2$, we have $\p^l_{n-1}\p^k_n=0$.
\end{theorem}
\begin{proof}
First we prove that $\mu_{n-1}()\mu_n\delta^k_n=\mu_{n-1}\p^k_n=0$, for any $n>1$ and $k\in\{1,\ldots,n-1\}$. It is a useful fact when defining the first homology group. One of the properties of the binary tree is that if the addresses in $T_n$, with $n>1$, are divided into fours, then the middle two addresses in each four have the same number of ones. Moreover, for each $k\in\{1,\ldots,n-1\}$, $n>1$, there is at least one four such that the middle two addresses have exactly $k$ 1s. Let $i$, $i+1$ be the positions of such addresses, so $T_n[i]_1=T_n[i+1]_1=k$. Then $(i,i+1)\in S^k_{2^n}$. The map 
$\sigma\mapsto (i,i+1)\sigma$ gives a bijection between even and odd permutations of $S^k_{2^n}$. The entropic condition (used after applying
$\mu_{n-1}\mu_n$) ensures that for any $2^n$-tuple $w$, we have $\mu_{n-1}\mu_n(i,i+1)\sigma w=\mu_{n-1}\mu_n\sigma w$ for any $\sigma\in S^k_{2^n}$. Thus,
$$\mu_{n-1}\mu_n\sum_{\sigma \in S^k_{2^n}} \sgn(\sigma) \sigma w=0,$$ that is, $\mu_{n-1}()\mu_n\delta^k_n w=0$.

Now we will show that for $k\in\{1,\ldots,n-1\}$, $l\in\{1,\ldots,n-2\}$, $n>2$, and any $\rho\in S^l_{2^{n-1}}$, we have $\mu_{n-1}\rho\mu_n\delta^k_n=\mu_{n-1}\rho\p^k_n=0$. In the infinite binary tree, every node with an address containing $l$ 1s has two branches leading from it; the left branch leads to a node with an address having $l$ 1s, the right branch leads to a node with $l+1$ 1s. The map $\mu_n$ corresponds to going up the binary tree, from $T_n$ to $T_{n-1}$. It follows that if $x_i*x_{i+1}$ is an element in a $2^{n-1}$-tuple
$\mu_n w$ that has a position with address containing $l$ 1s, then $w[i]_1=l$ and $w[i+1]_1=l+1$. Thus, any permutation $\rho\in S^l_{2^{n-1}}$ determines two permutations: $\rho_1\in S^l_{2^{n}}$ and $\rho_2\in S^{l+1}_{2^{n}}$ such that 
$\rho\mu_n=\mu_n\rho_1\rho_2$. It now follows that $$\mu_{n-1}\rho\mu_n\delta^k_n w=\mu_{n-1}\mu_n\rho_1\rho_2\sum_{\sigma \in S^k_{2^n}} \sgn(\sigma) \sigma w=$$$$\mu_{n-1}\mu_n\sum_{\sigma \in S^k_{2^n}} \sgn(\sigma)\rho_1\rho_2\sigma w=0.$$ The last equality is true because if both $\rho_1$ and $\rho_2$ do not permute addresses with $k$ 1s, then $\rho_1\rho_2\sigma w=\sigma\rho_1\rho_2 w$, and we can use the first part of the proof applied to a $2^n$-tuple $\rho_1\rho_2 w$. If one of them, say $\rho_1$, permutes the addresses with $k$ 1s, then 
$$\sum_{\sigma \in S^k_{2^n}} \sgn(\sigma)\rho_1\rho_2\sigma w= \sgn(\rho_1)\sum_{\sigma \in S^k_{2^n}}\sgn(\rho_1\sigma)\rho_1\sigma \rho_2 w,$$ and again the equality holds because we take the sum over all permutations from $S^k_{2^n}$.
Now, since $\delta^l_{n-1}$ in $\p^l_{n-1}$ is just a sum of signed permutations $\rho\in S^l_{2^{n-1}}$, it follows that $\p^l_{n-1}\p^k_n=0$.
\end{proof} 

For a map $\delta^\nu_n$ as defined before, let $\p^\nu_n=\mu_n\delta^\nu_n$.

\begin{corollary}
For any $\alpha\in\Z^{n-1}$ and $\beta\in\Z^{n-2}$, we have $\p^\beta_{n-1}\p^\alpha_n=0$.
\end{corollary}
\eop
By selecting a sequence of differentials: $d_n=\p^\nu_n$ for $n>1$, $d_1=\mu_1$, and $d_i=0$ for $i\leq 0$,
we obtain a chain complex $\{C_n(X),d_n\}$. The (co)homology of this chain complex is called the entropic (co)homology.
For an abelian group $G$, define the chain complex $C_*(X; G)=C_*\otimes G$, with $d'=d\otimes id$.
The groups of cycles and boundaries are denoted respectively by $Ker(d_n)=Z_n(X; G)\subset C_n(X; G)$
and $Im(d_{n+1})=B_n(X; G)\subset C_n(X; G)$.
The $n$-th entropic homology group of the entropic magma $(X;*)$ with coefficient group $G$ is defined as $$H_n(X; G)=H_n(C_*(X; G))=Z_n(X; G)/B_n(X; G).$$

\begin{lemma}
Let $w=(x_1,\ldots,x_{2^n})\in C_n$ satisfy the condition: there are $i$ and $j$, $i\neq j$, such that $w[i]_1=k$, $w[j]_1=k$, for some $k\in\{1,\ldots,n-1\}$, and $x_i=x_j$.
Then $\delta^k_n(w)=0$ (and thus also $\p^k_n(w)=0$).
\end{lemma}
\begin{proof}
If the above condition holds, then for any $\sigma\in S^k_{2^n}$ we have $\sigma (i,j)w=\sigma w$.
It follows that: $$\delta^k_n(w)=\sum_{\sigma \in S^k_{2^n}} \sgn(\sigma) \sigma w=\sum_{\substack{\sigma \in S^k_{2^n}\\\sgn(\sigma)=1}} \sigma w -
\sum_{\substack{\sigma \in S^k_{2^n}\\\sgn(\sigma)=-1}} \sigma w=$$$$
\sum_{\substack{\sigma \in S^k_{2^n}\\\sgn(\sigma)=1}} \sigma w - \sum_{\substack{\sigma \in S^k_{2^n}\\\sgn(\sigma)=-1}} \sigma (i,j)w=0,$$
because the second sum is now over all even permutations $\sigma (i,j)\in S^k_{2^n}$.
\end{proof}
\begin{corollary}
Let $\nu=(a_1,\ldots,a_{n-1})\in\Z^{n-1}$ and $w=(x_1,\ldots,x_{2^n})\in C_n$. Then $\p^\nu_n(w)=0$ if $w$ satisfy the following condition:
for every $i\in\{1,\ldots,n-1\}$ with $a_i\neq 0$, there are $x_{r_i}$ and $x_{s_i}$, $r_i\neq s_i$, such that $w[r_i]_1=i$, $w[s_i]_1=i$, and $x_{r_i}=x_{s_i}$.
\end{corollary}
\eop
For example, any tuple $(x_1,a,x_3,b,a,x_6,b,x_8)$ is a cycle in $Z_3(X; \Z)$.

It is useful to have a homomorphism $\xi^k_n\colon \Z X^{2^n}\to \Z X^{2^{n}}$, with $k\in\{1,\ldots,n-1\}$, defined on $2^n$-tuples $w$ by:
$$\xi^k_n w=(s,s+1)w + w,$$ where $s$ is the first place in $w$ such that $s=4i+2$, for some $i$, and $w[s]_1=w[s+1]_1=k$. In other words, for $w$ divided into fours, $\xi^k_n$ finds the first four in which the two middle places have $k$ ones in their addresses, permutes the corresponding elements, and adds $w$, so that $\xi^k_n w$ is a fixed element of the permutation $(s,s+1)$.

\begin{lemma}
For any $w=(x_1,\ldots,x_{2^n})\in C_n$, we have $\delta^k_n\xi^k_n w=0$ (and $\p^k_n\xi^k_n w=0$). Further,
for any $\nu=(a_1,\ldots,a_{n-1})\in\Z^{n-1}$, $\p^\nu_n\xi^{n-1}_n\ldots\xi^2_n\xi^1_n w=0$; here $\xi^i_n$ can be omitted if $a_i=0$.
\end{lemma}
\begin{proof}
$\delta^k_n\xi^k_n w=0$ because for any $\sigma\in S^k_{2^n}$, $\sigma(s,s+1)\xi^k_n w=\sigma\xi^k_n w$, and $\sigma\mapsto\sigma(s,s+1)$ is a bijection between even and odd permutations. This holds for any $k\in\{1,\ldots,n-1\}$ if we take $\xi^{n-1}_n\ldots\xi^2_n\xi^1_n w$ instead of $\xi^k_n w$, thus, the second statement is true.
\end{proof}
It follows that every chain from $C_n(X)$ can by made into a cycle by adjusting it with respect to every map $\delta^k_n$ appearing in the differential
$\p^\nu_n$. In particular, a sequence of elements of an entropic magma assigned to the last level of the resolving tree of a link (as described in Section
\ref{Entropic magma}) determines a cycle in the entropic homology of the magma. The topological meaning of this correspondence will be investigated in the future work, for now we consider an example. 

\begin{example}
Let ${\mathcal A}=(A,*,a_1,a_2,\ldots)$ be an entropic magma with a sequence of elements satisfying the conditions of Section \ref{Entropic magma}.
The last level of the resolving tree of the framed link \parbox{0.9cm}{\psfig{figure=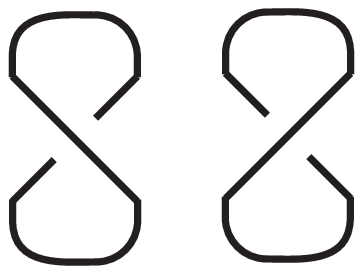,height=0.5cm}}, with the left crossing being resolved first, is: \parbox{2.5cm}{\psfig{figure=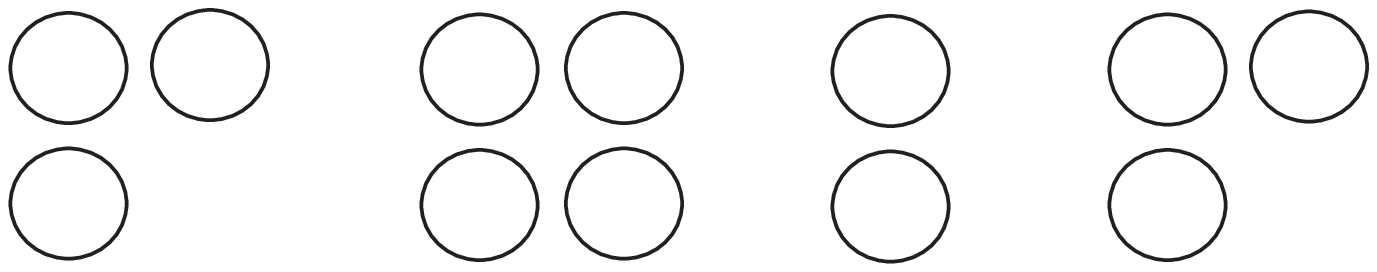,height=0.5cm}}. It corresponds to a list $w=(a_3,a_4,a_2,a_3)\in C_2(A)$.
In general, this list is not a cycle, but $\xi^1_2 w=(a_3,a_2,a_4,a_3)+(a_3,a_4,a_2,a_3)$ is:
$$\p^{(1)}_2\xi^1_2 w=\mu_2\delta^1_2\xi^1_2 w=\mu_2(()-(2,3))((2,3)w+w))=$$$$\mu_2(()(2,3)w+()w-(2,3)(2,3)w-(2,3)w)=\mu_2(0)=0.$$
Let ${\mathcal A}$ be the entropic magma

\begin{center}
\begin{tabular}{c| c cccc} 
$*$&1 & 2 & 3 & 4  \\
\hline 
1 & 2 & 1 & 3 & 4  \\
2 & 1 & 4 & 3 & 2  \\
3 & 3 & 3 & 3 & 3  \\
4 & 4 & 2 & 3 & 1  \\

\end{tabular}
\end{center}

\noindent with the sequence of elements $a_1=1$, $a_2=2$, $a_3=4$, $a_4=1$, $a_5=2$, $a_6=4, \ldots$. Then $w=(4,1,2,4)$ and calculations in GAP \cite{GAP} show that $\xi^1_2 w=(4,1,2,4)+(4,2,1,4)$ is one of the free generators of $H_2({\mathcal A})=Ker(\p^{(1)}_2)/Im(\p^{(1,0)}_3)$.

\end{example}

Now suppose that we are given a finite family $\{X_i=(X,*_i)\}$ of entropic magmas with an underlying set $X$, satisfying the compatibility condition:
$$(a*_ib)*_j(c*_id)=(a*_jc)*_i(b*_jd)$$ for any elements $a$, $b$, $c$, $d\in X$. In this situation we change the notation for the map $\mu_n$.
Let $\mu^{*_i}_n\colon \Z X^{2^n}\to \Z X^{2^{n-1}}$ be defined by:
$$\mu^{*_i}_n(x_1,x_2,x_3,x_4,\ldots,x_{2^n-1},x_{2^n})=(x_1*_ix_2,x_3*_ix_4,\ldots,x_{2^n-1}*_ix_{2^n}).$$
For $X_1,\ldots,X_r$ as above, let $\tau=(\tau_1,\ldots,\tau_r)\in\{+,-\}^r$ (basically, we are fixing the sign for each operation $*_i$).
Then we define:
$$\mu^\tau_n=\sum_{i\in\{1,\ldots,r\}}\tau_i\mu^{*_i}_n.$$
Also, let $\p^{k,*_i}_n=\mu^{*_i}_n\delta^k_n$ and $\p^{k,\tau}_n=\mu^\tau_n\delta^k_n$ for $k\in\{1,\ldots,n-1\}$, 
$\p^{\nu, *_i}_n=\mu^{*_i}_n\delta^\nu_n$ and $\p^{\nu,\tau}_n=\mu^{\tau}_n\delta^\nu_n$ for $\nu\in\Z^{n-1}$.
The map $\p^{\nu,\tau}_n$ is obtained by adding (or subtracting, depending on $\tau$) the considered earlier differentials $\p^\nu_n$ for each operation $*_i$, $i\in\{1,\ldots,r\}$. A simple check gives the following lemma.

\begin{lemma}
$\mu^\tau_1\p^{\nu,\tau}_2=\mu^\tau_1\mu^{\tau}_2\delta^\nu_2=0$.
\end{lemma}
\eop
This leads to the ``untwisted" (with $s=t=1$) version of the entropic cocycle condition that we considered at the beginning of this section:
$$f(x_1,x_2) - f(x_1,x_3) + f(x_3,x_4) - f(x_2,x_4) +$$$$
f(x_1*x_2,x_3*x_4) - f(x_1*x_3,x_2*x_4)=0.$$ Here we used three entropic operations: $*_1=*$, and the two trivial ones, $a*_Lb=a$, $a*_Rb=b$, for all $a$ and $b$, with $\tau=\{+,-,-\}$.

In general, $\mu^\tau_{n-1}\p^{\nu,\tau}_n\neq 0$. We can however make a change that leads to a setting in which this kind of condition holds.
Let us define a homomorphism $\zeta_n\colon \Z X^{2^n}\to \Z X^{2^{n}}$ defined on $2^n$-tuples by:
$$\zeta_n(x_1,x_2,x_3,x_4,\ldots,x_{4i},x_{4i+1},x_{4i+2},x_{4i+3},\ldots,x_{2^n-3},x_{2^n-2},x_{2^n-1},x_{2^n})=$$$$
(x_1,x_3,x_2,x_4,\ldots,x_{4i},x_{4i+2},x_{4i+1},x_{4i+3},\ldots,x_{2^n-3},x_{2^n-1},x_{2^n-2},x_{2^n}).$$
The involution $\zeta_n$ appears naturally when we have a family of compatible entropic magmas, namely:
$$\mu^{*_i}_{n-1}\mu^{*_j}_n w=\mu^{*_j}_{n-1}\mu^{*_i}_n \zeta_n w.$$
Now we change the maps $\delta^k_n(w)$. Let $$\hat{\delta}^k_n w =\sum_{\sigma \in S^k_{2^n}} \sgn(\sigma) (\sigma w - \zeta_n\sigma w).$$
We note that $\zeta_n\hat{\delta}^k_n w=-\hat{\delta}^k_n w$.

Because of the above change in the definition of $\delta^k_n(w)$, we slightly change the notation for the remaining maps:\\
for $k\in\{1,\ldots,n-1\}$: $\hat{\p}^{k,*_i}_n=\mu^{*_i}_n\hat{\delta}^k_n$, and $\hat{\p}^{k,\tau}_n=\mu^\tau_n\hat{\delta}^k_n$;\\
for $\nu=(a_1,\ldots,a_{n-1})\in\Z^{n-1}$:
$\hat{\delta}^\nu_n=\sum_{i\in\{1,\ldots,n-1\}}a_i\hat{\delta}^i_n$, 
$\hat{\p}^{\nu, *_i}_n=\mu^{*_i}_n\hat{\delta}^\nu_n$, and $\hat{\p}^{\nu,\tau}_n=\mu^{\tau}_n\hat{\delta}^\nu_n$.

\begin{theorem}
For any $\nu=(a_1,\ldots,a_{n-1})\in\Z^{n-1}$ and $\tau=(\tau_1,\ldots,\tau_r)\in\{+,-\}^r$, we have 
$\mu^\tau_{n-1}\hat{\p}^{\nu,\tau}_n=\mu^\tau_{n-1}\mu^{\tau}_n\hat{\delta}^\nu_n=0.$
\end{theorem}
\begin{proof}
For any $k\in\{1,\ldots,n-1\}$,
$$\mu^\tau_{n-1}\hat{\p}^{k,\tau}_n=\mu^\tau_{n-1}\mu^{\tau}_n\hat{\delta}^k_n=
\sum_{i\in\{1,\ldots,r\}}\tau_i\mu^{*_i}_{n-1}\sum_{i\in\{1,\ldots,r\}}\tau_i\mu^{*_i}_n\hat{\delta}^k_n.$$
We have: $$\mu^{*_i}_{n-1}\mu^{*_i}_n\hat{\delta}^k_n=\mu^{*_i}_{n-1}\mu^{*_i}_n\zeta_n\hat{\delta}^k_n=-\mu^{*_i}_{n-1}\mu^{*_i}_n\hat{\delta}^k_n.$$
Thus, $\mu^{*_i}_{n-1}\mu^{*_i}_n\hat{\delta}^k_n=0$.
Similarly, $$\mu^{*_i}_{n-1}\mu^{*_j}_n\hat{\delta}^k_n=\mu^{*_j}_{n-1}\mu^{*_i}_n\zeta_n\hat{\delta}^k_n=-\mu^{*_j}_{n-1}\mu^{*_i}_n\hat{\delta}^k_n.$$
Thus, $\mu^\tau_{n-1}\hat{\p}^{k,\tau}_n=0$. Also, $\mu^\tau_{n-1}\hat{\p}^{\nu,\tau}_n=\sum_{k\in\{1,\ldots,n-1\}}a_k\mu^\tau_{n-1}\mu^{\tau}_n\hat{\delta}^k_n$, so the result follows.
\end{proof}

Fixing a $\tau=(\tau_1,\ldots,\tau_r)\in\{+,-\}^r$, and choosing a sequence of homomorphisms $\hat{\p}^{\nu,\tau}_n$ for
$n>1$, allows us to define groups:
$$\hat{H}_n(X)=Ker(\mu^{\tau}_n)/Im(\hat{\p}^{\nu,\tau}_{n+1})=Ker(\mu^{\tau}_n)/Im(\mu^{\tau}_{n+1}\hat{\delta}^\nu_{n+1}),$$
for $n>0$, and 
$$\hat{H}_0(X)=C_0(X)/Im(\mu^{\tau}_1).$$
The maps involved in this definition are depicted in the following diagram.
\[\xymatrix @C=2.7pc{       & C_0\ar[ld]_0        & C_1\ar[ld]_{\mu^\tau_1} & C_2\ar[ld]_{\mu^\tau_2} & C_3\ar[ld]_{\mu^\tau_3} & \ar[ld] \ldots \\
 C_{-1} & C_0\ar[l]^{\mu^\tau_0=0}\ar[u]^(0.4){id} & C_1\ar[l]^{\mu^\tau_1}\ar[u]^(0.4){id} & C_2\ar[l]^{\mu^\tau_2}\ar[u]^(0.4){\hat{\delta}^\nu_2} & C_3\ar[l]^{\mu^\tau_3}\ar[u]^(0.4){\hat{\delta}^\nu_3} & \ar[l] \ldots }\]
We remark that we could also use $\tau'=(-1)^n\tau$ instead of $\tau$ in the above definition.

The calculations of $\hat{H}_n$ in GAP show that frequently these groups have torsion elements. We were able to calculate some groups $\hat{H}_1$ and $\hat{H}_2$, and the torsion parts included powers of $Z_2$, $Z_3$, and $Z_4$; we expect different $Z_k$'s for higher groups $\hat{H}_n$.

\section{Acknowledgements}
J.~H.~Przytycki was partially supported by the NSA-AMS 091111 grant and by the GWU REF grant.


\begin{thebibliography}{99}

\bibitem[C-Eil]{C-Eil} H.~Cartan, S.~Eilenberg, {\it Homological Algebra}, Princeton
University Press, 1956.

\bibitem[C-E-S]{C-E-S} J.~S.~Carter, M.~Elhamdadi, M.~Saito,
{\it Twisted quandle homology theory and cocycle knot invariants},
Algebraic \& Geometric Topology \textbf{2} (2002), 95-135.

\bibitem[Eth]{Eth}
M.~H.~Etherington, {\it Non-associative arithmetics}, Proc. Roy. Soc. Edinburgh \textbf{62A} (1949), 442-453.

\bibitem[GAP4]{GAP} 
The GAP Group, GAP -- Groups, Algorithms, and Programming,\\ 
{\tt  http://www.gap-system.org}

\bibitem[Lis]{Lis} J.~B.~Listing,  {\it  Vorstudien zur Topologie},  G\"ottinger Studien
(Abtheilung 1) \textbf{1} (1847),  811-875.

\bibitem [Mur]{Mur} D.~C.~Murdoch,
{\it Quasi-groups which satisfy certain generalized laws}, American
Journal of Math. \textbf{16} (1939), 509-522.

\bibitem [Prz-1]{Prz-1}
J.~H.~Przytycki, {\it Survey on recent invariants in classical
knot theory}, Warsaw University, Preprints 6,8,9, Warsaw, 1986.\\
e-print: \ {\tt http://front.math.ucdavis.edu/0810.4191}

\bibitem[Prz-2]{Prz-18}
J.~H.~Przytycki,
{\it Teoria w\c ez\l\'ow: podej\'scie kombinatoryczne},
(Knots: combinatorial approach to the knot theory),
Script, Warsaw, August 1995, 240+XLVIIIpp.

\bibitem[Prz-3]{Prz-37}
J.~H.~Przytycki,
{\it From Goeritz matrices to quasi-alternating links},
A Book Chapter (Chapter 9) in: The Mathematics of Knots: Theory and
Applications, Eds. M.~Banagl, D.~Vogel,
Springer-Verlag, 2011, 257-315.\\
e-print: \  {\tt  http://front.math.ucdavis.edu/0909.1118}

\bibitem[Prz-4]{Prz-100}
J.~H.~Przytycki, 
{\it {\bf KNOTS:} From combinatorics of knot diagrams to the
combinatorial topology based on knots}, Cambridge University Press,
accepted for publication, to appear 2013, pp. 600.\\
Chapter II, e-print:\ {\tt http://arxiv.org/abs/math/0703096}\\
Chapter III, e-print:\ {\tt arXiv:1209.1592v1 [math.GT]}\\
Chapter IV, e-print:\ {\tt http://front.math.ucdavis.edu/0909.1118}\\
Chapter V, e-print:\ {\tt http://arxiv.org/abs/math.GT/0601227}\\
Chapter VI, e-print:\ {\tt http://front.math.ucdavis.edu/1105.2238}\\
Chapter IX, e-print:\ {\tt http://arxiv.org/abs/math.GT/0602264}\\
Chapter X, e-print:\ {\tt http://arxiv.org/abs/math.GT/0512630}.

\bibitem [P-T-1]{P-T-1}
J.~H.~Przytycki, P.~Traczyk, {\it Invariants of links of Conway type},
Kobe J. Math. \textbf{4} (1987), 115-139.

\bibitem [P-T-2]{P-T-2}
J.~H.~Przytycki, P.~Traczyk, {\it Conway algebras and skein equivalence
of links}, Proc. Amer. Math. Soc. \textbf{100} (1987), no. 4, 744-748.

\bibitem[Sus]{Sus} A.~K.~Sushkevich, The theory of generalized groups, GNTI, Kharkov-Kiev, 1937 (In Russian).

\bibitem [Toy-1]{Toy-1}
K.~Toyoda, {\it On axioms of mean transformations and automorphic transformations of Abelian groups},
Tohoku Math. J. \textbf{46} (1940), 239-251.

\bibitem [Toy-2]{Toy-2} K.~Toyoda,
{\it On affine geometry of abelian groups}, Proceedings of the Imperial Academy \textbf{16} (1940), no. 5, 161-164.

\bibitem [Toy-3]{Toy-3} K.~Toyoda,
{\it On linear functions of abelian groups}, Proceedings of the Imperial Academy \textbf{16} (1940), no. 10, 524-528.

\bibitem [Toy-4]{Toy-4} K.~Toyoda,
{\it On axioms of linear functions}, Proceedings of the Imperial Academy \textbf{17} (1941), no. 7, 221-227.

\end{thebibliography}
\end{document}